\documentclass[11pt]{article}

\usepackage{amsfonts}
\usepackage{graphics}
\usepackage{amssymb}

\usepackage[dvips]{graphicx}

\newtheorem{theorem}{Theorem}[section]

\newcommand{\proof}{{\noindent \bf Proof:} }

\newcommand{\eop }{ \hfill $\Box$ }

\newcommand{\R}{\mathbf R}

\begin{document}

\begin{center}

\vspace{1cm}

 {\Large {\bf Application of an averaging principle on foliated diffusions: 
topology of the leaves.
\\[3mm]}}

\end{center}

\vspace{0.3cm}

\begin{center}
{\large {Paulo R.
Ruffino}\footnote{e-mail:
ruffino@ime.unicamp.br. Partially supported by FAPESP nr. 12/18780-0, 
11/50151-0 and CNPq nr. 477861/2013-0. 
}}

\vspace{0.2cm}

\textit{Departamento de Matem\'{a}tica, Universidade Estadual de Campinas, \\
13.083-859- Campinas - SP, Brazil.}

\end{center}

\begin{abstract}

We consider an $\epsilon K$ transversal perturbing vector field in a 
foliated Brownian motion defined in a foliated tubular neighbourhood of an 
embedded 
compact submanifold in $\R^3$. We study the effective behaviour of the system 
under this $\epsilon$ perturbation. If the perturbing vector field $K$ is 
proportional to the Gaussian curvature at the corresponding leaf, we have that 
the transversal component, after rescaling the time by $t/\epsilon$, approaches 
a linear increasing behaviour proportional to the Euler characteristic of $M$, 
as $\epsilon$ goes to zero. An estimate of the rate of convergence is presented.
\end{abstract}


\vspace{0.3cm}
\noindent {\bf MSC2010 subject classification:} 60H10, 58J65, 58J37.

\section{The set up}

The purpose of these notes is to explore a topological application of an 
averaging principle for foliated stochastic flows, as established in 
Gargate and Ruffino \cite{Gargate-Ruffino}. Our geometrical 
setting is a 
foliation of a 
tubular neighbourhood  of an embedded compact boundaryless submanifold $M$ 
of $\R^3$ consisting on diffeomorphic copies of $M$.
Precisely, consider a tubular neighbourhood $U$ of $M$ and a 
diffeomorphism 
$\phi:U \rightarrow M\times (-a, a) $ such that each $s\in 
(-a,a)\in \R$ 
generates a leaf $\phi^ {-1}(M \times \{s\} )$ which is 
diffeomorphic to $M= \phi^{-1} (M \times \{0\} )$. The 
diffeomorphism 
$\phi$ is taken such that transversal component $\frac{\partial 
\phi^{-1}}{\partial s}$ is orthogonal to the leaves, pointing outwards.

The main idea of this application is to 
consider an 
unperturbed system whose trajectories stay, each one, in a unique leaf and are 
Brownian motions on its corresponding leaf. This structure consisting of 
simultaneous 
Brownian motions on each leaf is called a {\it foliated Brownian motion}, see 
e.g. the seminal article by L. Garnett \cite{Garnett} or more recently 
Catuogno, Ledesma and Ruffino \cite{CRL-Trans}, \cite{CRL-Embedded}. We destroy 
this foliated behaviour of 
trajectories introducing   a
perturbing vector field orthogonal to the leaves such that at each point, this 
vector field
is given by the Gaussian curvature of 
the corresponding leaf, pointing outwards for positive curvature. Putting 
together the results on the averaging principle, \cite{Gargate-Ruffino}, 
for this particular foliated system and the classical Gauss-Bonnet theorem, we 
have 
that, in the 
average, the transversal behaviour of the perturbed system, with time 
rescaled by $\frac{t}{\epsilon}$, approaches a deterministic ODE with constant 
coefficient proportional to the Euler characteristic of the original 
submanifold $M$, 
as $\epsilon$ goes to zero. An estimate of the rate of convergence can be 
obtained.

In the next paragraphs we recall some of the results we are going to use 
on averaging stochastic flows in foliated spaces. 

\bigskip

\noindent{\it Previous results on averaging in foliated spaces}

\bigskip

 We recall the main results in \cite{Gargate-Ruffino} which are relevant here. 
Let $N$ be a smooth 
Riemannian manifold with an
$n$-dimensional
smooth foliation,
i.e. $N$ is endowed with an
integrable regular distribution of dimension $n$, hence $N$ is decomposed in a 
disjoint union of immersed submanifolds. For a precise definition
and further properties of foliated spaces in more general settings see e.g. the 
initial chapters of
Tondeur \cite{Tondeur}, Walcak
\cite{Walcak}, among others. We denote by $L_x$ the leaf 
of the foliation passing through a point $x\in N$.
We assume that the leaves are compact
and that given an initial condition $x_0$, the leaf  $L_{x_0}$ has a
tubular neighbourhood $U\subset N$ such that there exists a diffeomorphism 
$\varphi: U 
\rightarrow L_{x_0}\times V$, where $V\subset \R^d$ is a connected open set 
containing the origin. For simplicity, the second (vertical) coordinate of a 
point 
$x \in U$ is called the vertical projection $p(x)\in V$, i.e. 
$\varphi(x)=(u, \pi(x))$ for some $u\in L_{x}$.  
Hence for any fixed $v\in V$, the inverse image $p^{-1}(v)$ is a compact 
leaf $L_x$, where $x$ is any point in $U$ such that the 
vertical projection $p(x)=v$. In coordinates, we denote
\[
 p(x) = \Big(p_1(x), \ldots , p_d(x) \Big) \in V \subset \R^d
\]
for any $x\in U$. 

Consider an SDE in $M$ whose associated stochastic flow preserves the 
foliation, i.e.
we consider a Stratonovich equation

\begin{equation} \label{eq_original}
 dx_t = X_0 (x_t) dt + \sum_{k=1}^r X_k (x_t) \circ dB^k_t
\end{equation}
 where the smooth vector fields $X_k$ are foliated in the
sense that  $X_k(x) \in T_x L_x$,
for $k=0,1, \ldots , r$. Here $B_t = (B^1_t, \ldots , B^r_t)$ is a standard
Brownian motion in $\R ^r$
with respect to a filtered probability space $(\Omega, \mathcal{F}_t,
\mathcal{F}, \mathbf{P})$. For an 
initial condition $x_0$, the trajectories of the solution $x_t$ in this case
lay on the leaf $L_{x_0}$ a.s.. Moreover, there exists a (local) stochastic
flow of diffeomorphisms $F_t: N \rightarrow N$ which restricted
to the initial leaf is a flow in the compact submanifold $L_{x_0}$.

We introduce a perturbing smooth vector field $K$ in the system such that  
this vector field destroys 
the foliated structure of the trajectories. We denote the perturbed system by
$x^{\varepsilon}_t$ which satisfies the SDE
\begin{equation} \label{eq_perturbed}
 dx^{\varepsilon}_t = X_0 (x^{\varepsilon}_t) dt + \sum_{k=1}^r X_k
(x^{\varepsilon}_t) \circ dB^k_t + \varepsilon K (x^{\varepsilon}_t)\, dt,
\end{equation}
with the same initial condition $x^{\varepsilon}_0=x_0$.

\bigskip

For each vertical coordinate $i=1, \ldots, d$, we denote by $Q^i(v)$ the 
ergodic average of the $i$-th component of the perturbing vector field $K$ in 
the corresponding leaf $p^{-1}(v)$. Hence, by the ergodic theorem
\[
 Q^i (v) := \int_{p^{-1}(v)} Dp_i (K)(x)\ d\mu_v (x) = \lim_{t \rightarrow 
\infty} \frac{1}{t}\int_0^t Dp_i (K)(F_s (x))\ ds
 \]
for $\mu_v$-almost every point $x $ in the leaf $p^{-1}(v)$; here $D_xp_i$ 
is the derivative of $p_i$ at the point $x$ and $\mu_v$ denotes an invariant 
measure on $p^{-1}(v)$ which we assume uniquely ergodic. Denote by 
$\eta(t)$ an estimate of the rate of convergence of the ergodic limit in $L^p$ 
for the functions $Q^i(v)$. Hence $\eta(t)$ tends to zero when $t$ goes to 
infinity. In general, there is no optimal rate of convergence, see e.g. 
Kakutani-Petersen \cite{Kakutani-Petersen},
Krengel \cite{Krengel} and an explicit example of continuous system in the 
averaging context in 
\cite{Gargate-Ruffino}.

We assume that the averaging functions of the perturbing vector field  $Q^i: 
V \rightarrow \R$ are Lipschitz continuous for $i=1, \ldots, d$.
This hypothesis holds naturally if the invariant measures $\mu_v$ for the
unperturbed foliated system has a sort of weakly continuity on $v$. For 
example, for nondegenerate systems on the leaves, this condition is naturally 
satisfied.

Consider the following ODE on the vertical space $V$:
\[
 \frac{ \mathrm{d} v(t)}{\mathrm{d} t} = \Big( Q^1(v(t)), \ldots , Q^d(v(t)) 
\Big),
\]
with initial condition $p(x_0)=0$. Let $T_0$ be the time that the 
solution $v(t)$ hits the boundary $\partial V$. Let 
$\tau^\epsilon$ be the stopping time given by the exit time of the perturbed 
system $x^{\varepsilon}_t$ from the coordinate neighbourhood $U \subset M$.

An averaging principle in this context is established by the theorem 
below which says that the transversal 
behaviour of
$x^{\varepsilon}_{\frac{t}{\epsilon}}$  can be approximated in the average by an
ordinary differential equation in the transversal space  whose
coefficients are given by the average of the transversal component of the
perturbation $K$ with respect to the invariant measure on the leaves for the
original dynamics of equation (\ref{eq_original}), when $\epsilon$ decreases to 
zero.
The rate of converge is given below:

%
%

\begin{theorem} \label{teoremaprincipal} Assuming Lipschitz continuity of the 
averaging functions  $Q^i: 
V \rightarrow \R$  for $i=1, \ldots, d$
we have:

\begin{description}
 \item[(1)] For any  $0<t<T_0$, $\epsilon>0$, $\beta \in (0,1/2)$, $\alpha \in 
(0,1)$ and  $2
\leq q< \infty$, there exist 
functions $C_1=C_1(t)$ and $C_2=C_2(t)$
such that
\[   
 \left[\mathbb{E}\left( \sup_{s\leq t}\left|p
\left( x^{\varepsilon}_{\left( \frac{s}{\epsilon}\right) \wedge \tau^{\epsilon}}
\right)-v (s)\right|^q\right)\right]^{\frac{1}{q}} 
 \leq  C_ 1 \, \epsilon^{\alpha} +  C_2\,  \eta \left( t |\ln 
\epsilon |^{\frac{2\beta}{q}} \right),
\]
 where $\eta (t)$ is the rate of convergence in $L^q$ of the 
ergodic averages 
of the unperturbed trajectories on the leaves.

\item[(2)] For $\gamma > 0 $, let 
\[
 T_{\gamma} = \inf \  \{t>0 \ |\ \mathrm{dist} (v(t), \partial V ) \leq 
\gamma \}.
\]
The exit times of the two systems satisfy the estimates

\[
\mathbb{P}(\epsilon \tau^{\epsilon}<T_{\gamma})\leq
\gamma^{-q} \left[  C_ 1(T_{\gamma}) \, \epsilon^{\alpha} +  C_2 (T_{\gamma})\, 
\eta \left( T_\gamma |\ln 
\epsilon |^{\frac{2\beta}{q}} \right)\right]^q.
\]

\end{description}

\end{theorem}

Item (b) of the theorem above guarantees robustness of the result. For 
the proof of the theorem above see \cite{Gargate-Ruffino}; further extension to 
L\'evy 
processes has been done in \cite{Hoegele-Ruffino}.

\section{Exploring the topology of the leaves}


We consider an unperturbed foliated 
dynamics in $U$ (degenerate) which is a 
foliated Brownian motion. That is, for each initial condition 
$x_0\in U$, the solution is a Brownian motion on the corresponding leaf 
$L_{x_0}$. For more details and construction of this processes, see 
\cite{CRL-Trans}. In our particular 
case of embedded manifold we can consider the dynamics generated by 
gradient vector fields tangent to the leaves. Precisely, at each  $x\in U$, let 
$X^i(x)$ be the orthogonal projection of $e_i$, the $i$-th element of the 
canonical basis onto the tangent space $T_x L_{x}$ of the leaf passing through 
$x$, for each $i=1,2,3$. The vector fields $X^1, X^2 $ and $X^3$ determine the 
following Stratonovich stochastic equation in $U$:
\[
 dx_t = \sum_{i=1}^3 X^i(x_t) \circ dB_t^i.
\]
The corresponding stochastic flow of this equation generates a 
foliated Brownian motion in $U$, 
i.e. given  initial conditions, the solution are simultaneous Brownian motions
on each leaf of the foliation, see \cite{CRL-Embedded}.

We investigate the effective behaviour of a small
transversal perturbation of order $\epsilon$:
\[
 dx_t^\epsilon = \sum_{i=1}^3 X^i(x_t^\epsilon) \circ dB_t^i + \epsilon 
K(x_t^\epsilon),
\]
where the transversal vector field $K(x)$ is orthogonal to the leaves, 
proportional to the Gaussian curvature of the leaf $L_x$ at the point $x\in U$. 
We assume that at a point with positive curvature, the corresponding vector 
field $K$ points outwards.  In some sense 
this equation models an average inertial or centrifugal forces acting on 
particles 
moving randomly on the leaves.
We have the following averaging result which does not depend on the geometry of 
the 
manifold $M$:

\begin{theorem} \label{Thm: Euler}  The transversal dynamics 
of orthogonal perturbation, given by  
$\epsilon$ times the Gaussian curvature of foliated 
Brownian motion, behaves according 
to 
$v(t)= 2 \pi \chi(M)\ 
t$, as $\epsilon$ goes to zero. Here 
$\chi(M)$ is the Euler characteristic of the leaves.

Precisely, up to a 
stopping time, for $q\geq 2$, $\beta \in (0,1/2)$, $\epsilon>0$, for each 
$t\geq 0$, there exists a constant $C>0$ such that 
\begin{equation} \label{eq:rate of convergence}
 \left[\mathbb{E}\left( \sup_{s\leq t}\left|p
\left( x^{\epsilon}_{\left( \frac{s}{\epsilon}\right)}
\right)-v (s)\right|^q\right)\right]^{\frac{1}{q}} < C |\ln 
\epsilon |^{-\frac{\beta}{q}}.
\end{equation}

\end{theorem}

\proof Most of the proof is straightforward: one just has to put the 
properties together and apply the classical Gauss-Bonnet theorem. In fact, the 
invariant measures on the leaves are the normalized Riemannian 
volume. According to our coordinate systems, the one dimensional 
transversal 
component is given by the Gaussian curvature $K(x)$. The average (in the 
transversal direction) of the 
perturbing vector field $K$ on each leaf is given by
\[
 \int_{L_x} K(x) \ d\mu (x) =2 \pi \chi(M).
\]
Hence this average depends only on the topology of the leaf. 

The proof finishes 
as a direct application of Theorem 
\ref{teoremaprincipal} with codimension $d=1$.  The rate of ergodic 
convergence on the leaves $\eta(t)$, in this case, has order $1/\sqrt{t}$ 
since the 
system is uniformly elliptic, as in X.-M.-Li \cite{Li}. The estimate on 
the rate of convergence (inequality (\ref{eq:rate of convergence})) follows 
directly  .

\eop

 As a simple example, 
consider an orthogonal perturbation, according to the curvature, of a Brownian 
motion on 
any 
manifold diffeomorphic to the torus. This system presents no transversal 
behaviour on 
the 
average. Indeed, a flat torus illustrates trivially this fact since $K\equiv 
0$. 
On 
manifolds diffeomorphic to the 
sphere, the  orthogonal behaviour is a linear expansion. For an
$n$-fold torus with genus $n>1$, the 
effect is a linear contraction with coefficients 
proportional to $2-2n$.

\end{document}